\documentclass[10pt]{amsart}
\usepackage{color,amsmath,amssymb,graphics,amsthm}
\usepackage{hyperref}
\usepackage{enumerate}
\newtheorem{theorem}{Theorem}

\newtheorem{conjecture}[theorem]{Conjecture}

\newtheorem{definition}[theorem]{Definition}

\newtheorem{example}[theorem]{Example}

\newtheorem{lemma}[theorem]{Lemma}

\def\ie{{\em i.e. }}

\def\>{{\rangle}}
\def\<{{\langle}}

\newcommand{\FF}{{\mathbb{F}}}
\newcommand{\Q}{\mathbb{Q}}

\newcommand{\Ext}{\mathop{\mathrm{Ext}}\nolimits}
\newcommand{\Irr}{\mathop{\mathrm{Irr}}\nolimits}

\newcommand{\Hecke}{{\mathcal{H}}}

\newcommand{\blambda}{\boldsymbol{\lambda}}

\include{abbr}

\date{}
\newcommand{\cnode}[1]{\kern1pt\rlap{\hbox{$\kern 3pt\scriptstyle #1$}}
\kern-4pt\mathop\bigcirc\kern-3.5pt}

\title[Decomposition matrices for $d$-Harish-Chandra series]{Decomposition matrices for $d$-Harish-Chandra series: the exceptional rank two cases}
\author{Maria Chlouveraki}
\address[M. Chlouveraki]{
	School of Mathematics, University of Edinburgh, JCMB, King's Buildings,
	Edinburgh, EH9 3JZ, UK
}
\email{maria.chlouveraki@ed.ac.uk}
\author{Hyohe Miyachi}
\address[H. Miyachi]{
Graduate School of Mathematics, Nagoya University, Chikusa-ku, Nagoya
Aichi, 464-8602, Japan
}
\email{miyachi@math.nagoya-u.ac.jp}

\begin{document}
\maketitle

\begin{abstract}
We  calculate all decomposition matrices of the cyclotomic Hecke
 algebras of the rank $2$ exceptional complex reflection groups in
 characteristic $0$. We prove the existence of canonical basic sets in
 the sense of Geck-Rouquier and show that all modular irreducible representations
 can be lifted to the ordinary ones.
\end{abstract}

\section{Introduction}
\subsection{Notation}

Throughout this paper $W$ will be a complex reflection group and $K$ will be its field of definition, \ie the smallest field to which the traces of all the elements of $W$ belong.  We denote by $\mu(K)$ the group of all roots of unity in the number field $K$.
For any positive integer $d$, we set
$\zeta_d:=\mathrm{exp}(2\pi i/d)$ and we denote by  $\Phi_d(x)$ the $d$-th cyclotomic polynomial over $\Q$. From now on, let $d$ be a positive integer and let $q$ be an indeterminate.

\subsection{Motivation from $d$-Harish-Chandra theory}\label{mdHC}

In \cite{BMR},
the \emph{cyclotomic Hecke algebra} $\Hecke_q(W)$ associated with the complex reflection group $W$  is introduced: it is defined over $\mathbb{C}[q,q^{-1}]$ as a quotient algebra of the group algebra of the braid group associated with $W$.
The motivation for introducing this algebra 
 is that
$\Hecke_q(W)$ should appear naturally 
as an endomorphism algebra of
$R_L^G(\blambda)$,
where $R_L^G$ is the Deligne-Lusztig twisted induction
 functor for the finite reductive group $G$ and the so-called $d$-cuspidal 
pair $(L,\blambda)$ (cf.~\cite{BMM1}).
The complex $R_L^G(\blambda)$ is a trunction of the one giving a
so-called Brou\'e abelian defect conjecture.
In particular, if we take $d=1$ and $\blambda=1$, this picture goes back to
 the classical one, where
 $R_L^G(\blambda)$ is the permutation module over the (finite) Borel
 subgroup of $G$ and the Hecke algebra is  the
 classical Iwahori-Hecke algebra.
 Hence, $\Hecke_q(W)$ should be treated as the endomorphism algebra
 of the $d$-analogue of the permutation module over the Borel (principal
 series modules), \ie the $d$\emph{-Harish-Chandra series.}

The important points for $\Hecke_q(W)$ are that
\begin{itemize}
 \item 
 this algebra is a flat deformation of the group algebra of $W$
 (actually, for $q=\zeta_d$, $\Hecke_q(W)$ specializes to the group algebra $\mathbb{C}W$),
 \item 
 the relations deforming the group algebra multiplication (``cyclotomic relations'')
could involve multiparameters, 
 \item
 there should be some distinguished choices of those parameters
 which are matched up with the endomorphism ring of the 
complex that we have seen above,
 \item those special choices for parameters are not known in general.
{\footnote{
The authors are indebted to C. Bonnaf\'e for informing them of this 
open problem.
}}
\end{itemize}

\subsection{Aim and results}\label{tableau}

In this paper, we shall consider the {\it{asymptotic}} modular representation theory of
the $d$-Harish-Chandra series 
in the case where
the {\it{conjectural}} choices for the distinguished parameters are
 known and $W$ is an exceptional (irreducible) complex reflection group of rank $2$.
Here, by ``asymptotic'' we mean the
 modular representation theory in characteristic zero.
 The conjectural choices of the parameters by Brou\'e and Malle (cf.~\cite{BrMa},  \cite{Malle}) are given in the following table:  \footnote{Some cases are omitted, since they can obtained from others by replacing $q$ with $-q$. Moreover, we 
 include the cases corresponding to $H_4$ for the same reasons that they are included in \cite{BrMa} and  \cite{Malle}.}

\begin{center}
$ $\\
\begin{tabular}{|c|c|c|c|c|}
  \hline
  $W$ & $G$ & $d$ & $\blambda$ & $\textrm{Parameters of } \Hecke_q(W)$ \\ 
  \hline
  $G_{4}$ & $^3D_4$  & $3$ & $1$ & $1,q,q^2$ \\
  \hline
   $G_{5}$ & $F_4$  & $3$ & $1$ & $1,q,q^2\,;\,1,q,q^2$ \\
    & $^2E_6$  & $3$ & $1$ & $1,q,q^2\,;\,1,q^2,q^4$ \\
   & $E_8$  & $3$ & $^3D_4[-1]$ & $1,q,q^2\,;\,1,q^4,q^8$ \\
  \hline
    $G_{8}$ & $^2F_4$  & $\,\,8'$ & $1$ & $1,\zeta_8^3q,\zeta_8^5q,q^2$ \\
    & $F_4$  & $4$ & $1$ & $1,q,-q,q^2$ \\
   & $^2E_6$  & $4$ & $1$ & $1,q,q^2,q^3$ \\
    & $E_7$  & $4$ & $\phi_2^3$ & $1,q,-q,-q^4$ \\
    &               & $4$ & $\phi_{11}^3$ & $1,q^3,-q^3,-q^4$ \\
    & $E_8$  & $4$ & $\phi_{3,1}$ & $1,-1,-q,q^5$ \\
    &               & $4$ & $\phi_{123,013}$ & $1,-q^4,q^5,-q^5$ \\
    &               & $4$ & $\phi_{013,2}$ & $1,-q,-q^4,q^5$ \\
\hline
    $G_{9}$ & $E_8$  & $8$ & $1$ & $1,q^4\,;\,1,q^2,q^4,q^6$ \\
 \hline
    $G_{10}$ & $E_8$  & $12$ & $1$ & $1,-q^2,q^4\,;\,1,q^3,-q^3,q^6$ \\
  \hline
    $G_{12}$ & $^2F_4$  & $4$ & $1$ & $1,q^2$ \\
    \hline
    $G_{16}$ & $E_8$  & $5$ & $1$ & $1,q,q^2,q^3,q^4$ \\
    \hline
    $G_{20}$ & $H_4$  & $3$ & $1$ & $1,q,q^2$ \\
    \hline
    $G_{22}$ & $H_4$  & $4$ & $1$ & $1,q^2$ \\         
\hline   
\end{tabular}
$ $\\$ $
\end{center}

We will calculate the decomposition matrices for all of the above cases and show the existence of \emph{canonical basic sets}, as defined by Geck and Rouquier. The existence of canonical basic sets has already been proved for all Weyl groups for all choices of parameters in characteristic $0$ (cf.~\cite{GR}, \cite{Gsurvey}, \cite{GJ}, \cite{ChlJa}). In that sense, $d$-Harish-Chandra theory works as the ordinary Harish-Chandra theory.
 What is more, we will see that, for the exceptional complex reflection groups of rank $2$, all ``modular'' irreducible representations are obtained as modular reductions of some ``ordinary'' ones, which therefore comprise an \emph{optimal basic set}.

\section{Preliminaries}

\subsection{Cyclotomic Hecke algebras} 
Let $V$ be a finite dimensional $K$-vector space such that $W \subset \mathrm{GL}(V)$ acts irreducibly on $V$. 
Let us denote by $\mathcal{A}$ the set of reflecting hyperplanes of $W$ in $V$ and set
$\mathcal{M}:=\mathbb{C} \otimes V - \bigcup_{H \in \mathcal{A}} \mathbb{C} \otimes H$. For $x_0 \in \mathcal{M}$, let $B=\Pi_1(\mathcal{M}/W,x_0)$ be the braid group associated with $W$.

For every orbit $\mathcal{C}$ of $W$ on $\mathcal{A}$, we set
$e_{\mathcal{C}}$ to be the common order of the subgroups $W_H$, where $H$
is any element of $\mathcal{C}$ and $W_H$ the subgroup formed by $\mathrm{id}_V$
and all the reflections fixing the hyperplane $H$.

We choose a set of indeterminates
$\textbf{u}=(u_{\mathcal{C},j})_{(\mathcal{C} \in
\mathcal{A}/W)(0\leq j \leq e_{\mathcal{C}}-1)}$ and we denote by
$\mathbb{Z}[\textbf{u},\textbf{u}^{-1}]$ the Laurent polynomial ring
in all the indeterminates $\textbf{u}$. We define the \emph{generic
Hecke algebra} $\mathcal{H}(W)$ of $W$ to be the quotient of the group
algebra $\mathbb{Z}[\textbf{u},\textbf{u}^{-1}]B$ by the ideal
generated by the elements of the form
$$(\textbf{s}-u_{\mathcal{C},0})(\textbf{s}-u_{\mathcal{C},1}) \ldots (\textbf{s}-u_{\mathcal{C},e_{\mathcal{C}}-1}),$$
where $\mathcal{C}$ runs over the set $\mathcal{A}/W$ and
$\textbf{s}$ runs over the set of monodromy generators around the
images in $\mathcal{M}/W$ of the elements of the hyperplane
orbit $\mathcal{C}$.\footnote{The monodromy generators are the generators of the braid group $B$. For their definition, the reader may refer, for example, to \cite[\S 4.1.2]{springer}.}

Let $\zeta_d \in \mu(K)$. A $\zeta_d$-\emph{cyclotomic specialization} is a morphism of $\mathbb{C}$-algebras of the form
$$\begin{array}{rccc}
\varphi: & \mathbb{C}[\textbf{u},\textbf{u}^{-1}] & \rightarrow & \mathbb{C}[q,q^{-1}]\\ \smallbreak \smallbreak
& u_{\mathcal{C},j} & \mapsto & \zeta_{e_\mathcal{C}}^j (\zeta_d^{-1}q)^{n_{\mathcal{C},j} }
\end{array}$$
with $n_{\mathcal{C},j}  \in \mathbb{Z}$ for all $\mathcal{C},j$. The algebra $\mathcal{H}_\varphi(W)$ obtained as the specialization of $\mathcal{H}(W)$ via the morphism $\varphi$ is a $\zeta_d$-\emph{cyclotomic Hecke algebra}. For $q = \zeta_d$,  $\mathcal{H}_\varphi(W)$ specializes to the group algebra $\mathbb{C}W$. A \emph{cyclotomic Hecke algebra} is a Hecke algebra which is $\zeta_d$-cyclotomic for some positive integer $d$.

Now, in the spirit of generalized Harish-Chandra theory, $\zeta_d$-cyclotomic Hecke algebras may correspond to $d$-Harish-Chandra series for some finite reductive group. 
The table of Section \ref{tableau} gives us a list of all 
cyclotomic Hecke algebras associated with exceptional complex reflection groups of rank $2$ which should appear as $d$-Harish-Chandra series for some finite reductive group $G$. Each one of these Hecke algebras is $\zeta_d$-cyclotomic for the corresponding $d$.

\subsection{Schur elements and generic degrees}
\def\sfD{\mathsf{D}}

From now on, $\mathcal{H}_q$ will denote a cyclotomic Hecke algebra associated with the group $W$.
Set $y^{|\mu(K)|}:=q$. The algebra 
$\mathbb{C}(y) \mathcal{H}_q:=
\mathbb{C}(y)
\otimes_{{\mathbb{C}}[q,q^{-1}]}\mathcal{H}_q$
 is split semisimple (\cite[Proposition 4.3.4]{springer}) and by Tits's deformation theorem, there exists a bijection 
$\mathrm{Irr}(\mathbb{C}(y)\mathcal{H}_q) \rightarrow \mathrm{Irr}(W),\,\chi_q \mapsto \chi$. 

Furthermore, it is conjectured
that $\mathcal{H}_q$ has a symmetrizing form $t_q: \mathcal{H}_q \rightarrow \mathbb{C}[q,q^{-1}]$ (that is, a trace function such that
the bilinear form 
$(h,h')\mapsto t_q(hh')$ is non-degenerate) with nice properties (cf.~\cite[Theorem 2.1]{BMM2}): for example, $t_q$ specializes to the canonical symmetrizing form on $\mathbb{C}W$ for $q=\zeta_d$. This conjecture has been verified for all but a finite number of cases (see \cite{MaMi} for the list of cases). In particular, it has been verified for all groups studied in this paper, except for $G_{16}$, $G_{20}$ and $G_{22}$.

Due to Geck \cite{Ge}, we know that
$$t_q=\sum_{\chi \in \mathrm{Irr}(W)}\frac{1}{s_{\chi}}\chi_q,$$
where $s_{\chi}$ is the \emph{Schur element} of $\mathcal{H}_q$
associated with
 $\chi \in \mathrm{Irr}(W)$. The Schur elements belong to $\mathbb{Z}_K[y,y^{-1}]$ and they are products of cyclotomic polynomials over $K$ (cf.~\cite[Proposition 4.3.5]{springer}). For all $ \chi \in \mathrm{Irr}(W)$, we set
$$a_\chi\,:=\frac{\mathrm{valuation}(s_\chi(y))}{|\mu(K)|}.$$

Finally, let $P(y) \in K[y,y^{-1}]$ be a product of cyclotomic polynomials over $K$ which is divisible by all Schur elements of $\Hecke_q$. Following ~\cite[\S 6E]{BMM2}, we
define the $P$-\emph{generic degree} of $\mathcal{H}_q$ associated with $\chi \in \mathrm{Irr}(W)$ to be the quotient
$$D^{P}_{\chi}:=\frac{P(y)}{s_\chi}.$$
If we take $P(y)$ to be the \emph{relative Poicar\'e polynomial} of $\mathcal{H}_q$ (in the sense of Brou\'e-Malle-Michel ~\cite[\S 6A]{BMM2}), then $D^{P}_{\chi}$ is called simply the \emph{generic degree}. If $\blambda=1$ in the corresponding cuspidal pair, then the relative Poicar\'e polynomial is equal to the Schur element associated with the trivial character (up to multiplication by a power of $y$).

\subsection{Decomposition matrix and basic sets}

Let $\theta:\mathbb{C}[y,y^{-1}]\to \mathbb{C}$ be a ring homomorphism 
such that
$\theta (q )=\xi \in \mathbb{C}^\times$. Considering $\mathbb{C}$ as a $\mathbb{C}[y,y^{-1}]$-module via $\theta$, we set
$\mathcal{H}_{\xi}:=\mathbb{C} \otimes_{\mathbb{C}[y,y^{-1}] }\mathcal{H}_{q}$. We call $\theta$ a \emph{specialization} of $\mathcal{H}_q$.
Note that the algebra $\mathcal{H}_\xi$ is split.

   Recall that there is a bijection $\mathrm{Irr}(\mathbb{C}(y)\mathcal{H}_q) \leftrightarrow \mathrm{Irr}(W)$.    Let 
$R_0(\mathbb{C}(y)\mathcal{H}_{q})$ (respectively $R_0(\mathcal{H}_{\xi})$) denote the Grothendieck group
 of finitely generated $\mathbb{C}(y)\mathcal{H}_{q}$-modules (respectively $\mathcal{H}_{\xi}$-modules). It is generated by the classes $[U]$ of the simple  $\mathbb{C}(y)\mathcal{H}_{q}$-modules (respectively $\mathcal{H}_{\xi}$-modules) $U$.
Following \cite[Theorem 7.4.3]{GePf}, we obtain a well-defined  decomposition map
$$d^{\theta}:R_0 (\mathbb{C}(y)\mathcal{H}_{q})    \to R_0 (\mathcal{H}_{\xi})$$
with corresponding decomposition matrix $(d_{\chi,\phi}^\theta)_{\chi \in \mathrm{Irr}(W),\phi \in \mathrm{Irr}(\Hecke_\xi)}$. 
The decomposition matrix records information about the irreducible representations of the specialized Hecke algebra $\mathcal{H}_\xi$ in terms of the
irreducible representations of the split semisimple Hecke algebra $\mathbb{C}(y)\mathcal{H}_{q}$. In particular, it records how the irreducible representations of $\mathbb{C}(y)\mathcal{H}_{q}$ break into irreducible representations of $\mathcal{H}_\xi$.

If $\mathcal{H}_\xi$ is semisimple, then Tits's deformation theorem implies that there exists a canonical bijection between the set of irreducible characters of $\mathbb{C}(y)\mathcal{H}_q$ and the set of irreducible characters of $\mathcal{H}_\xi$. Hence, the corresponding decomposition matrix is the identity matrix and the simple modules of $\mathcal{H}_\xi$ are parametrized by $\mathrm{Irr}(W)$. However, what do we do when the algebra $\mathcal{H}_\xi$ is not semisimple?

 The theory of basic sets gives standard ways to parametrize the simple modules for Hecke algebras. 
 A \emph{basic set} is a subset of $\mathrm{Irr}(W)$ which is in bijection with the set of irreducible characters of $\mathcal{H}_\xi$. We will be mostly interested in two kinds of basic sets: the ``canonical'' basic sets (in the sense of Geck-Rouquier \cite{GR}) and the ``optimal'' ones.
 
\begin{definition}\label{canbasic}
\emph{We say that $\mathcal{H}_q$  admits a \emph{canonical basic set} $\mathcal{B}^{\mathrm{can}} \subset \mathrm{Irr}(W)$ with respect to $\theta$  if and only if the following two conditions are satisfied:
  \begin{enumerate}[(1)]
\item   For all $\phi \in\operatorname{Irr}(\mathcal{H}_\xi)$, there exists $\chi_\phi\in \mathcal{B}^{\mathrm{can}}$ such that
\begin{enumerate}[(a)]
 \item $d^\theta_{\chi_\phi,\phi} =1$;
 \item $\textrm{if }d^\theta_{\psi,\phi} \neq 0 \textrm{ for some $\psi \in \mathrm{Irr}(W)$, then either }
 \psi=\chi_\phi \textrm{ or } a_\psi>a_{\chi_\phi}$. \footnote{The original definition of the canonical basic set uses Lusztig's $a$-function, which is defined as the opposite of the $a$-function that we use here. In the case of classical Iwahori-Hecke algebras of Weyl groups, the original definition yields that the trivial character has minimal $a$-value in $\mathrm{Irr}(W)$ (which is equal to $0$). However, this fact depends on the order that we write down the parameters of our algebra. Therefore, in the cases that will be studied in this paper, since we copy the order of parameters from \cite{BrMa} and \cite{Malle}, we prefer to use the $a$-function defined in Section 2.1, so that the trivial character has always minimal $a$-value. This choice does not affect the existence of canonical basic sets (see also \cite{ChlJa}), but may affect the parametrization of the characters belonging to them.}
 \end{enumerate}
 \item The map $\mathrm{Irr}(\mathcal{H}_\xi)  \mapsto \mathcal{B}^{\mathrm{can}},\,\phi \mapsto \chi_\phi$ is a bijection.
\end{enumerate}}
\end{definition}

Hence, if $\mathcal{H}_q$  admits a canonical basic set, then we can reorder the rows of the decomposition matrix so that its upper part (consisting of the rows indexed by the elements of
$\mathcal{B}^{\mathrm{can}}$) is a unitriangular matrix.
 
 The existence of canonical basic sets in characteristic $0$ has been proved for all choices of parameters for Weyl groups (\cite{GR}, \cite{Gsurvey}, \cite{GJ}, \cite{ChlJa}), for the Ariki-Koike algebras (\cite{DJM}, \cite{GR}, \cite{Ariki}, \cite{Ug}, \cite{Jacon}) and for certain cyclotomic Hecke algebras associated with the groups $G(de,e,r)$ (\cite{GeJa}).

\begin{definition}\label{optimalbasic}
\emph{We say that $\mathcal{H}_q$  admits an \emph{optimal basic set} $\mathcal{B}^{\mathrm{opt}} \subset \mathrm{Irr}(W)$ with respect to $\theta$  if and only if the following two conditions are satisfied:
  \begin{enumerate}[(1)]
\item   For all $\phi \in\operatorname{Irr}(\mathcal{H}_\xi)$, there exists $\chi_\phi\in \mathcal{B}^{\mathrm{opt}}$ such that
 \begin{enumerate}[(a)]
 \item $d^\theta_{\chi_\phi,\phi} =1$;
 \item $\textrm{if }d^\theta_{\psi,\phi} \neq 0 \textrm{ for some $\psi \in \mathrm{Irr}(W)$, then either }
 \psi=\chi_\phi \textrm{ or } \psi \notin  \mathcal{B}^{\mathrm{opt}}.$
 \end{enumerate}
 \item The map $\mathrm{Irr}(\mathcal{H}_\xi)  \mapsto \mathcal{B}^{\mathrm{opt}},\,\phi \mapsto \chi_\phi$ is a bijection.
\end{enumerate}}
\end{definition}

Hence, if $\mathcal{H}_q$  admits an optimal basic set, then all irreducible representations of ${\mathcal{H}}_\xi$ are obtained as modular reductions of irreducible representations of $\mathcal{H}_q$: we say that all irreducible representations of $\Hecke_\xi$ can be \emph{lifted} to $\Hecke_q$.
Therefore, we can reorder the rows of the decomposition matrix so that its upper part (consisting of the rows indexed by the elements of
$\mathcal{B}^{\mathrm{opt}}$) is the identity matrix.
Moreover, the optimal basic set produces the full decomposition matrix automatically.

In all the cases studied in this paper, we will show that
$\mathcal{H}_q$ admits an optimal basic set with respect to any specialization $\theta$
(if $\mathcal{H}_\xi$ is semisimple, there is nothing to prove).
Moreover, we have checked (using GAP) that $\mathcal{H}_q$ also
admits a canonical basic set with respect to any $\theta$. 

\subsection{Block partitions and good central elements}\label{good central}
Let $\theta:\mathbb{C}[y,y^{-1}]\to \mathbb{C},\,q \mapsto \xi$ be a specialization of $\mathcal{H}_q$ as above. We say that
$\chi,\, \psi \in \mathrm{Irr}(W)$ \emph{belong to the same block} if they label the rows of the same block in the decomposition matrix $$(d_{\chi,\phi}^\theta)_{\chi \in \mathrm{Irr}(W),\phi \in \mathrm{Irr}(\Hecke_\xi)}.$$ If an irreducible character is alone in its block, then we call it a {\em character of defect zero}. We have that  $\chi$ is an irreducible character of defect zero if and only if
 $\theta(s_\chi) \neq 0$ (see \cite[Lemma 2.6]{MaRo}).

 It is well known that
 $\chi,\psi \in \Irr(W)$ belong to the same block if and only if
 $\theta(\omega_\chi(z)) = \theta(\omega_\psi(z))$
 for any $z \in Z(\Hecke_q)$, where $\omega_\chi,\, \omega_\psi$ denote the corresponding central characters (see \cite[Lemma 7.5.10]{GePf}).
Unfortunately, any precise description of the centre, as in \cite{GeRo},
 is not yet known for our cyclotomic Hecke algebras.
However, by \cite[4.20, 5.13]{MR1429870},
we can know the character values of certain ``good'' central elements.
Here is a list of some good central elements for our complex reflection groups:
\[
\begin{array}{c|ccccccc}
{\mbox{Type}}&G_4&G_5&G_8&G_9&G_{10}&G_{16}&G_{20}\\
{\mbox{Good element}}&(st)^3&(st)^2&(st)^3 &(st)^3&(st)^2&(st)^3&(st)^5\\
\end{array} 
\]

Let $\chi \in \Irr(W)$. Following the results of Brou\'e and Michel,
 the value of $\omega_\chi$ at these good central elements is of the form
\begin{equation*}
\lambda(\chi):=
{\mathrm{exp}}\left(
2\pi \sqrt{-1} (a_\chi+A_{\chi}) /d^2)\right)  \gamma_{\chi} q^{(a_\chi+A_\chi)/d},
\end{equation*}
where  $d$ and $\gamma_\chi$ are complex numbers completely determined by $W$.
So  we get a necessary condition
\begin{equation*}\theta(\lambda(\chi))=\theta(\lambda(\psi)) 
\end{equation*}
 for the block linkage of $\chi,\, \psi \in \mathrm{Irr}(W)$.  Note that this argument is not new:  Geck used this method in \cite[12.4]{MR1179271} for the case of Weyl groups.

\begin{example} 
Let $W=G_{10}$ and
$$\Hecke_q(G_{10})=\left\langle s,t\,\left|\, stst=tsts, 
\begin{array}{l}
(s-1)(s+q^2)(s-q^4)=0,\\ (t-1)(t- q^3)(t+  q^3)(t-q^6)=0\end{array} \right.\right\rangle.$$
Following the notation of characters in the \emph{GAP3} package {\em CHEVIE}, each irreducible character of $W$ is denoted by $\phi_{d,b}$, where $d$ is the dimension of the representation and $b$ is the valuation of its fake degree.
The value $\lambda(\phi_{d,b})$
 of $\Hecke(G_{10})$ at the central element 
$(st)^2=(ts)^2$
is given by the following table:
{\small{
\[
\begin{array}{c|c|cc|cc|cc|cccc}
%0&4&5&5&6&6&7&7&8&8&8&8\\
%0&46&57&57&68&68&78&78&84&88&88&88\\
\phi_{1,0}&\phi_{1,8}&\phi_{2,1}&\phi_{2,4}&\phi_{1,6}&\phi_{1,12}&\phi_{2,5}&
\phi_{2,8}&\phi_{2,7}'&\phi_{3,2}&\phi_{2,7}''&\phi_{3,6}'\\\hline
1&q ^4&q ^5&-q ^5&q ^6&q ^6&-q ^7&q ^7&q ^8&\zeta_3^2 q ^8&-q ^8&q ^8\\
\end{array} 
\]
\[
\begin{array}{cc|cccc|cccccc}
%8&8&9&9&9&9&10&10&10&10&10&10\\
%88&88&95&95&96&96&100&104&104&104&104&104\\
\phi_{1,16}&\phi_{3,10}''&\phi_{4,9}&\phi_{4,3}&\phi_{2,9}&\phi_{2,12}&
\phi_{2,11}''&\phi_{3,4}&\phi_{2,11}'&\phi_{3,10}'&\phi_{1,14}&\phi_{3,12}''\\\hline
q ^8&\zeta_3 q ^8&-q ^9&q ^9&q ^9&-q ^9&q ^{10}&\zeta_3 q ^{10}&
-q ^{10}& \zeta_3 q ^{10}&q ^{10}&q ^{10}\\
\end{array}
\]
\[
\begin{array}{cccccc|cc|cccc}
%10&10&10&10&10&10&11&11&11&11&12&12\\
%104&104&104&104&104&104&107&107&108&108&108&112\\
\phi_{1,20}&\phi_{4,11}&\phi_{3,14}&\phi_{4,5}&\phi_{3,6}''&\phi_{3,8}'&
\phi_{4,13}&\phi_{4,7}&\phi_{2,13}&\phi_{2,10}&\phi_{2,15}'&\phi_{3,8}''\\\hline 
q ^{10}&-i q ^{10}&\zeta_3^2
q ^{10}&i q ^{10}&q ^{10}&\zeta_3^2 q ^{10}&q ^{11}&
-q ^{11}&-q ^{11}&q ^{11}&q ^{12}&\zeta_3^2
q ^{12}\\
\end{array}
\]
\[
\begin{array}{cccc|cc|cc|cc|c|c}
%12&12&12&12&13&13&14&14&15&15&16&20\\
%112&112&112&112&114&114&116&116&117&117&118&120\\
\phi_{3,12}'&\phi_{3,16}&\phi_{1,18}&\phi_{2,15}''&\phi_{2,14}&\phi_{2,17}&
\phi_{1,22}&\phi_{1,28}&\phi_{2,18}&\phi_{2,21}&\phi_{1,26}&\phi_{1,34}\\\hline 
q ^{12}&\zeta_3 q ^{12}&q ^{12}&-q ^{12}&-q ^{13}&q ^{13}&
q ^{14}&q ^{14}&q ^{15}&-q ^{15}&q ^{16}&q^{20}\\
\end{array}
\]
}}

Let $\theta$ be a specialization of $\Hecke_q(G_{10})$ such that $\theta(q)=\zeta_7$. By looking at the Schur elements, we can easily check (with the use of {\em GAP}) that all irreducible characters are of defect $0$, except for the following ones:
$$\phi_{1,0},\,\,\phi_{1,12},\,\,\phi_{2,8},\,\,\phi_{2,17}, \,\,\phi_{1,28}, \,\,\phi_{1,34}.$$ 
According to the above table, we have
$$\theta(\lambda(\phi_{1,0}))=\theta(\lambda(\phi_{2,8}))=\theta(\lambda(\phi_{1,28}))=1$$
and
$$\theta(\lambda(\phi_{1,12}))=\theta(\lambda(\phi_{2,17}))=\theta(\lambda(\phi_{1,34}))=\zeta_7^6.$$
Due to the necessary condition for block linkage given above, we must have two separate blocks  consisting of three characters each, namely
$$B_1=\{\phi_{1,0},\, \phi_{2,8},\,\phi_{1,28}\} \,\,\text{ and }\,\,
B_2=\{\phi_{1,12},\, \phi_{2,17},\,\phi_{1,34}\}.$$ 

\end{example}

\section{Methodology and results}

\subsection{Summarized main results and sketch of the proof}\label{criteria}

Let $\Hecke_q(W)$ be one of the cyclotomic Hecke algebras in
 the list of \S \ref{res}.  Set $\mathcal{O}:=\mathbb{C}[q,q^{-1}]$ and let $\theta: \mathcal{O} \rightarrow \mathbb{C},\,q \mapsto \xi$ be a ring homomorphism as above.

\begin{theorem}\label{main01}
There exists an optimal basic set for  $\Hecke_q(W)$ with respect to any specialization $\theta$.
\end{theorem}

We have proved the above result case by case, by calculating the decomposition matrix in all the cases where the algebra $\Hecke_\xi$ is not semisimple, that is all the cases where the decomposition matrix is not the identity matrix. By \cite[Theorem 7.4.7]{GePf}, $\Hecke_\xi$ is not semisimple if and only if there exists $\chi \in \mathrm{Irr}(W)$ such that $\theta(s_\chi)=0$. Due to the form of the Schur elements, this happens only if
$\xi$ is a root of unity. In \S \ref{res}, we give the optimal basic sets in all the cases that need to be considered.

Now, the following criteria have been enough to help us determine the decomposition matrices in all the cases which are of interest to us:
\begin{enumerate}[(1)]
\item We have $\theta(s_\chi) \neq 0$ if and only if $\chi$ is an irreducible character of defect zero.
\item All $1$-dimensional representations are irreducible.
\item If a representation has one of the $1$-dimensional subrepresenations found in the previous step, then it is reducible. In particular, $2$-dimensional representations are irreducible, unless they have an $1$-dimensional subrepresentation. This can be easily checked, because $1$-dimensional subrepresentations correspond to $1$-dimensional eigenspaces for the generating matrices of the $2$-dimensional representation.
\item We have the following necessary condition for block linkage (see \S \ref{good central}): If  $\chi,\, \psi \in \mathrm{Irr}(W)$ are in the same block,
      then $\theta(\lambda(\chi))=\theta(\lambda(\psi))$. 
\item We apply Lemma \ref{GR}, taking $P(y)$ to be the least common multiple of all Schur elements of $\Hecke_q$. 
\end{enumerate}
The above criteria give us in each case an optimal basic set and thus a complete description of the irreducible representations of $\mathcal{H}_\xi$.
We complete the decomposition matrix by using the fact that the modular reductions of the irreducible characters of $\mathcal{H}_q$ can be written uniquely as $\mathbb{N}$-linear combinations of the
 irreducible characters of $\mathcal{H}_\xi$.

The irreducible representations of the cyclotomic Hecke algebras associated with the exceptional complex reflection groups have been calculated by Malle and Michel in \cite{MaMi}. We are grateful to them for providing us with preliminary versions of their paper and for implementing everything (generating matrices, character tables, etc.) in the GAP3 package CHEVIE, which has been used for all our calculations.

Now, having computed the decomposition matrix, it is easy to check the validity of the following result:

\begin{theorem}\label{main02}
There exists a canonical basic set for  $\Hecke_q(W)$ with respect to any specialization $\theta$.
\end{theorem}

In the Appendix of this paper, we will work out in detail the example of $G_{12}$ for $\xi=\zeta_8$.
It is a small, but representative example of the computational work we had to do to obtain Theorems \ref{main01} and \ref{main02}, as well as the data presented in Section \ref{res}. It illustrates the method, the CHEVIE commands and the programmes used in the process.

\subsection{Decomposition matrix and generic degrees}

Recall that $q=y^{|\mu(K)|}$.
 Let $\zeta$ be a primitive root of unity of order $m$ and consider a specialization $\theta$ such that $\theta(y)=\zeta$. Let $\mathfrak{n}$ be the ideal of $\mathcal{O}$ generated by the irreducible polynomial $y-\zeta$. Note that  
$\mathcal{O}_\mathfrak{n}$ is a discrete valuation ring.

If $\zeta'$ is a primitive root of unity of order $m$ such that $\zeta \neq \zeta'$, and $\mathfrak{n}'$ is the ideal of $\mathcal{O}$ generated by  $y-\zeta'$, then the rings $\mathcal{O}_\mathfrak{n}$ and
 $\mathcal{O}_{\mathfrak{n}'}$ are isomorphic (see, for example, \cite[Proposition 1.1.2]{springer}).
Since the decomposition map $d^\theta$ is actually defined as the decomposition map of the algebra $\mathcal{O}_\mathfrak{n}\mathcal{H}_q$ (\cite[Theorem 7.4.3]{GePf}), the decomposition matrix of $\mathcal{H}_q$ with respect to $\theta$ is the ``same''\footnote{The indexing of the rows may be different: characters of the same degree may be interchanged.} as the decomposition matrix with respect to $\theta': y \mapsto \zeta'$. Hence, the form of the decomposition matrix does not depend on the choice of the root of unity $\zeta$, but on its order $m$.

Let $(d_{\chi,\phi}^\theta)_{\chi \in \mathrm{Irr}(W),\phi \in \mathrm{Irr}({\mathcal{H}_\xi})}$ be the decomposition matrix of $\Hecke_q$ with respect to $\theta$.
By  \cite[Proposition 4.4]{GeRo}, for all $\phi \in \mathrm{Irr}({\mathcal{H}_\xi})$, we must have
$$\sum_{\chi \in \mathrm{Irr}(W)}
\frac{d_{\chi,\phi}^\theta}{s_\chi} \in \mathcal{O}_\mathfrak{n}.$$

If $P(y) \in K[y,y^{-1}]$ is a product of $K$-cyclotomic polynomials  which is divisible by all Schur elements of $\Hecke_q$, then,  for all $\phi \in \mathrm{Irr}({\mathcal{H}_\xi})$,
we obtain that
$$\sum_{\chi \in \mathrm{Irr}(W)}
\frac{d_{\chi,\phi}^\theta D_\chi^P}{P(y)} \in \mathcal{O}_\mathfrak{n}.$$
We deduce that the polynomial
$$\sum_{\chi \in \mathrm{Irr}(W)}
{d_{\chi,\phi}^\theta D_\chi^P}$$ has at least as many factors of the form $(y-\zeta)$ as $P(y)$.
In particular, the following holds:

\begin{lemma}\label{GR}
Let $\theta:q \mapsto \xi$ be a specialiazation of $\mathcal{H}_q$ such that ${\mathcal{H}}_\xi$ is not semisimple and let $P(y)$ be an element of $K[y,y^{-1}]$ as above.
If $\zeta^{|\mu(K)|}=\xi$, then, for all $\phi \in \mathrm{Irr}({\mathcal{H}_\xi})$,
$$\sum_{\chi \in \mathrm{Irr}(W)}
{d_{\chi,\phi}^\theta D_\chi^P}(\zeta)=0.$$
\end{lemma}

\subsection{Decomposition matrices for the exceptional groups of rank $2$}\label{res}

Let $W$ be an exceptional irreducible complex reflection group of rank $2$. For each cyclotomic specialization appearing in \cite{BrMa} and \cite{Malle},
we will present here an optimal basic set for $\mathcal{H}_q(W)$ with respect to any specialization
$q=\xi$ such that  $\mathcal{H}_\xi$ is not semisimple. 
The modular reductions of the elements of the optimal basic set are the irreducible characters of $\Hecke_\xi$.
The actual
decomposition matrix can be calculated very easily, 
since all modular reductions of irreducible characters of
$\mathcal{H}_q$
 are written uniquely as ${\mathbb{N}}$-linear combinations of the elements of
 $\mathrm{Irr}(\Hecke_\xi)$.

We will not write down the characters of the basic set whose Schur
element does not become zero when $q=\xi$, \ie the characters of defect $0$. In the cases where we have more
than one optimal basic set (different irreducible characters of $\Hecke_q$ have the same irreducible modular reduction), we will choose the representatives with  minimal
$a$-value. We will also separate the blocks using parentheses.

Now, if $m$ is an odd positive integer and $\xi$ is a primitive root of unity of order $m$, then $-\xi$ is a primitive root of unity of order $2m$. This is why, in the cases when $\mathcal{H}_q \cong \mathcal{H}_{-q}$ (namely the cases (ii), (iv), (v) and (vii) for $G_8$ and the cases of $G_9$, $G_{10}$, $G_{12}$  and $G_{22}$), we will only study one of the two specializations.

Finally, we want to record some information on block structures which appear very often. Therefore, we will write:
\begin{center}
(i) $[\phi_1]$,\,\,\, 
(ii) $[\phi_1,\phi_3]$,\,\,\,
(iii) $\langle \phi_1,\phi_3 \rangle$,\,\,\,
(iv) $\{\phi_1,\phi_2\}$\,,\,\,
(v) $[[\phi_1,\phi_4, \phi_6]]$
\end{center}
if the decomposition matrix of the corresponding block is, respectively, of the form:

\begin{center}

$ $\\
(i) $\begin{array}{c|cc}
\phi_1 &1&  \\
\phi_2 &1& \\
\end{array}$\,\,\,\,\, \,\,\,\,\, 
(ii)
$\begin{array}{c|cc}
\phi_1 &1&  \\
\phi_2 &1&1 \\
\phi_3 & &1\\
\end{array}$
\end{center}
$ $\\
\begin{center} 
(iii)
$\begin{array}{c|cc}
\phi_1 &1&  \\
\phi_2 &1&1 \\
\phi_3 & &1\\
\phi_4 &1 &1 \\
\phi_5 &1 & \\
\end{array}$ \,\,\,\,\,\,\,\,\,\, \,\,\,\,\, 
(iv)
$\begin{array}{c|cc}
\phi_1 &1&  \\
\phi_2 & &1 \\
\phi_3 &1&1\\
\phi_4 &  &1 \\
\phi_5 &1 & \\
\end{array}$ \,\,\,\,\, \,\,\,\,\, \,\,\,\,\, 
(v)
$\begin{array}{c|ccc}
\phi_1 &1& & \\
\phi_2 &1&1& \\
\phi_3 &1&1&1\\
\phi_4 & &1& \\
\phi_5 & &1 &1 \\
\phi_6 & & &1 \\
\end{array}.$ 
\end{center}

$ $

The vacant entry in the decomposition matrices above stands for zero.
Actually, in the first two cases, by \cite[Theorem 3.6]{Uno}, we know that the block algebras 
with such decomposition matrices  are Brauer tree algebras.

The notation for the characters is the one used by the GAP3 package
CHEVIE. The trivial character is always denoted by $\phi_{1,0}$.

\subsubsection{Optimal basic sets for the group $G_4$}
$ $\\
\\
Let
$$\Hecke_q(G_4)=\left\langle s,t\,\left|\, sts=tst, 
\begin{array}{l}
(s-1)(s-q)(s-q^2)=0,\\ (t-1)(t-q)(t-q^2)=0\end{array} \right.\right\rangle.$$
The characters which, together with all characters of defect $0$, form an optimal basic set for $\Hecke_q(G_4)$ are:

$q=1:$   ($\phi_{1,0}),   (\phi_{2,1}$).

$q=-1:$  \{$ \phi_{1,0},   \phi_{2,1}$\}.

$q=\zeta_6:$  [[$ \phi_{1,0},  \phi_{1,4}, \phi_{1,8}$]].

$q=\zeta_{12}:$  [$ \phi_{1,0},  \phi_{1,8}$].

\subsubsection{Optimal basic sets for the group $G_5$}
$ $\\
\\
(i) Let
$$\Hecke_q(G_5)=\left\langle s,t\,\left|\, stst=tsts, 
\begin{array}{l}
(s-1)(s-q)(s-q^2)=0,\\ (t-1)(t-q)(t-q^2)=0\end{array} \right.\right\rangle.$$
The characters which, together with all characters of defect $0$, form an optimal basic set for $\Hecke_q(G_5)$ are:

$q=1:$   $(\phi_{1,0}),\,   (\phi_{2,1})$.

$q=-1:$  ($ \phi_{1,0}, \phi_{1,4'}, \phi_{1,4''}$),  $(\phi_{2,1})$.

$q=i:$  $ [[\phi_{1,0},  \phi_{1,8''},\phi_{1,16}]],\,
 [\phi_{1,4'}, \phi_{1,12'}],\, 
[\phi_{1,4''}, \phi_{1,12''}]$.

$q=\zeta_{6}:$  $ [[\phi_{1,12'}, \phi_{1,0},\phi_{1,12''}]],\, [[\phi_{1,4'},  \phi_{1,16}, \phi_{1,4''}]]$.

$q=\zeta_{8}:$  $[\phi_{1,0},  \phi_{1,16}]$.

$q=\zeta_{12}:$  $[\phi_{1,0},  \phi_{2,9}],\, [\phi_{2,1}, \phi_{1,16}]$.\\
$ $\\
(ii) Let
$$\Hecke_q(G_5)=\left\langle s,t\,\left|\, stst=tsts, 
\begin{array}{l}
(s-1)(s-q)(s-q^2)=0,\\ (t-1)(t-q^2)(t-q^4)=0\end{array} \right.\right\rangle.$$
The characters which, together with all characters of defect $0$, form an optimal basic set for $\Hecke_q(G_5)$ are:

$q=1:$   $(\phi_{1,0}),\,   (\phi_{2,1})$.

$q=-1:$  $ (\phi_{1,0}, \phi_{1,4''}),\, ( \phi_{2,3'}) $.

$q=i:$  $ \{\phi_{1,0},\phi_{1,8'''}\},\,[\phi_{2,1}],\, \{\phi_{2,3'},\phi_{1,4''}\},\,
 [\phi_{2,5'}]$.

$q=\zeta_6:$  $[[\phi_{1,0},  \phi_{1,8''},\phi_{1,16}]],\,
 [[\phi_{1,4'}, \phi_{1,12'}, \phi_{1,8'''}]],\, 
[[ \phi_{1,8'},\phi_{1,4''},\phi_{1,12''}]]$.

$q=\zeta_8:$  $ [\phi_{1,0}, \phi_{1,12''}],\, [\phi_{1,4'}, \phi_{1,16}]$.

$q=\zeta_{10}:$  $ [\phi_{1,0}, \phi_{1,12'}],\, [\phi_{1,4''}, \phi_{1,16}]$.

$q=\zeta_{12}:$  $ [\phi_{1,0}, \phi_{1,16}],\, [\phi_{1,4''}, \phi_{2,7''}],\,
[\phi_{1,12'}, \phi_{2,3'}]$.

$q=\zeta_{18}:$  $[\phi_{1,0},  \phi_{2,9}],\, [\phi_{2,1}, \phi_{1,16}]$.\\
$ $\\
(iii) Let
$$\Hecke_q(G_5)=\left\langle s,t\,\left|\, stst=tsts, 
\begin{array}{l}
(s-1)(s-q)(s-q^2)=0,\\ (t-1)(t-q^4)(t-q^8)=0\end{array} \right.\right\rangle.$$
The characters which, together with all characters of defect $0$, form an optimal basic set for $\Hecke_q(G_5)$ are:

$q=1:$   $(\phi_{1,0}),\,   (\phi_{2,1})$.

$q=-1:$  $ (\phi_{1,0},   \phi_{1,4''}),\, (\phi_{2,3'})$.

$q=i:$ $ (\phi_{1,0}, \phi_{1,8'''}),\,  (\phi_{1,4''}),\,
(\phi_{2,1}),\, (\phi_{2,5'})$.

$q=\zeta_6:$ $[[ \phi_{1,0}, \phi_{1,12'}, \phi_{1,12''}]],\,  [[ \phi_{1,8'}, \phi_{1,8''},\phi_{1,8'''}]],\, 
[[\phi_{1,4'}, \phi_{1,4''}, \phi_{1,16}]]$.

$q=\zeta_8:$  $ [\phi_{1,0}],\, \{\phi_{1,4''}, \phi_{2,3'}\},\, [\phi_{2,1}],\,
[\phi_{1,8'''}],\, [\phi_{2,5'}]$.

$q=\zeta_{10}:$  $ [[\phi_{1,0},  \phi_{1,8''}, \phi_{1,16}]],\, [\phi_{1,4'}, \phi_{1,12'}],\,
[\phi_{1,4''}, \phi_{1,12''}]$.

$q=\zeta_{12}:$  $[\phi_{1,0},   \phi_{1,12''}], [\phi_{1,8'},
\phi_{1,8'''}],\, [\phi_{1,16}, \phi_{1,4'}]$.

$q=\zeta_{14}:$ $ [\phi_{1,8'}, \phi_{1,4''}],\,  [\phi_{1,12'},\phi_{1,8'''}]$.

$q=\zeta_{18}:$ $ [\phi_{1,0}, \phi_{1,12'}],\,  [\phi_{1,4''},
\phi_{1,16}],\, [\phi_{1,8'}, \phi_{2,5'}],\, [\phi_{1,8'''},\phi_{2,5'''}]$.

$q=\zeta_{20}:$  $[\phi_{1,0},  \phi_{1,16}]$.

$q=\zeta_{24}:$ $ [\phi_{2,7''}, \phi_{1,4''}],\,  [\phi_{1,12'},\phi_{2,3'}]$.

$q=\zeta_{30}:$  $[\phi_{1,0},  \phi_{2,9}],\, [\phi_{2,1}, \phi_{1,16}]$.

\subsubsection{Optimal basic sets for the group $G_8$}
$ $\\
\\
(i) Let
$$\Hecke_q(G_8)=\left\langle s,t\,\left|\, sts=tst, 
\begin{array}{l}
(s-1)(s- \zeta_8^3 q)(s- \zeta_8^5 q)(s-q^2)=0,\\ (t-1)(t- \zeta_8^3 q)(t- \zeta_8^5 q)(t-q^2)=0\end{array} \right.\right\rangle.$$
The characters which, together with all characters of defect $0$, form an optimal basic set for $\Hecke_q(G_8)$ are:

$q=1$  : $\{\phi_{1,0},\phi_{3,2}\},$
$[\phi_{2,1}],$
$[\phi_{2,4}]$.

$q=-1$  : $\{\phi_{1,0},\phi_{3,2}\},$
$[\phi_{2,1}],$
$[\phi_{2,4}].$

$q=i$  : $[[\phi_{1,0},\phi_{2,7''},\phi_{1,18}]].$

$q=\zeta_{8}^{3}$  : $\langle \phi_{1,0},\phi_{2,10} \rangle,$
$\langle \phi_{1,6},\phi_{2,4} \rangle,$
$(\phi_{2,1}).$

$q=\zeta_{8}^{5}$  : $\langle \phi_{1,0},\phi_{2,13} \rangle,$
$\langle \phi_{1,12},\phi_{2,1} \rangle,$
$(\phi_{2,4}).$

$q=\zeta_{12}$  : $[\phi_{1,0},\phi_{1,18}].$

$q=\zeta_{24}$  : $[\phi_{1,0},\phi_{2,10}],$
$[\phi_{2,1},\phi_{2,13}],$
$[\phi_{2,4},\phi_{1,18}].$\\
$ $\\
(ii) Let
$$\Hecke_q(G_8)=\left\langle s,t\,\left|\, sts=tst, 
\begin{array}{l}
(s-1)(s+q)(s- q)(s-q^2)=0,\\ (t-1)(t+ q)(t-  q)(t-q^2)=0\end{array} \right.\right\rangle.$$
The characters which, together with all characters of defect $0$, form an optimal basic set for $\Hecke_q(G_8)$ are:

$q=-1$  : $(\phi_{1,0},\phi_{2,4}),$
$(\phi_{2,1}).$

$q=\zeta_{6}$  : $[[\phi_{1,0},\phi_{1,12},\phi_{1,18}]].$

$q=\zeta_{8}$  : $[\phi_{1,0},\phi_{3,8}],$
$[\phi_{3,2},\phi_{1,18}].$

$q=\zeta_{12}$  : $[\phi_{1,0},\phi_{1,18}],$
$[\phi_{2,1},\phi_{2,13}],$
$[\phi_{2,4},\phi_{2,10}].$\\
$ $\\
(iii) Let
$$\Hecke_q(G_8)=\left\langle s,t\,\left|\, sts=tst, 
\begin{array}{l}
(s-1)(s-q)(s- q^2)(s-q^3)=0,\\ (t-1)(t- q)(t-  q^2)(t-q^3)=0\end{array} \right.\right\rangle.$$
The characters which, together with all characters of defect $0$, form an optimal basic set for $\Hecke_q(G_8)$ are:

$q=1$  : $(\phi_{1,0},\phi_{3,2}),$
$(\phi_{2,1}).$

$q=-1$  : $(\phi_{1,0},\phi_{1,6},\phi_{2,1}).$

$q=\zeta_{3}$  : $\{\phi_{2,1},\phi_{2,4}\},$
$\{\phi_{1,0},\phi_{3,2}\}.$

$q=\zeta_{6}$  : $(\phi_{1,0},\phi_{1,6},\phi_{1,12},\phi_{1,18}),$
$[\phi_{2,4},\phi_{2,10}].$

$q=\zeta_{8}$  : $[\phi_{1,0},\phi_{2,10}],$
$[\phi_{2,4},\phi_{1,18}].$

$q=\zeta_{10}$  : $[\phi_{1,0},\phi_{2,13}],$
$[\phi_{2,1},\phi_{1,18}].$

$q=\zeta_{12}$  : $[\phi_{1,0},\phi_{1,12}],$
$[\phi_{1,6},\phi_{1,18}],$
$[\phi_{2,1},\phi_{2,13}].$

$q=\zeta_{18}$  : $[\phi_{1,0},\phi_{1,18}].$\\
$ $\\
(iv) Let
$$\Hecke_q(G_8)=\left\langle s,t\,\left|\, sts=tst, 
\begin{array}{l}
(s-1)(s-q)(s+q)(s+q^4)=0,\\ (t-1)(t- q)(t+q)(t+q^4)=0\end{array} \right.\right\rangle.$$
Surprisingly, the Poincar\'e polynomial is equal to the Schur element associated with $\phi_{1,18}$.
The characters which, together with all characters of defect $0$, form an optimal basic set for $\Hecke_q(G_8)$ are:

$q=-1$  : $(\phi_{1,0},\phi_{1,6},\phi_{2,1}).$

$q=\zeta_{3}$  : $(\phi_{1,0},\phi_{1,12},\phi_{2,7''}),$
$[\phi_{2,1},\phi_{2,13}],$
$[\phi_{3,2}].$

$q=\zeta_{5}$  : $[\phi_{1,0},\phi_{2,10}],$
$[\phi_{3,2},\phi_{1,18}].$

$q=\zeta_{7}$  : $[\phi_{2,4},\phi_{1,18}].$

$q=\zeta_{8}$  : $[\phi_{1,0}],$
$[\phi_{2,1}],$
$[\phi_{2,4}],$
$[\phi_{3,2}].$

$q=\zeta_{9}$  : $[\phi_{1,6},\phi_{1,18}].$

$q=\zeta_{12}$  : $[\phi_{1,0},\phi_{1,18}],$
$[\phi_{2,1},\phi_{2,13}],$
$[\phi_{2,4},\phi_{2,10}].$\\
$ $\\
(v) Let
$$\Hecke_q(G_8)=\left\langle s,t\,\left|\, sts=tst, 
\begin{array}{l}
(s-1)(s-q^3)(s+q^3)(s+q^4)=0,\\ (t-1)(t- q^3)(t+  q^3)(t+q^4)=0\end{array} \right.\right\rangle.$$
Surprisingly  in this case also, the Poincar\'e polynomial is equal to $q^{63}s_{\phi_{1,0}}$. The characters which, together with all characters of defect $0$, form an optimal basic set for $\Hecke_q(G_8)$ are:

$q=-1$  : $(\phi_{1,0},\phi_{1,6},\phi_{2,1}).$

$q=\zeta_{3}$  : $(\phi_{1,0},\phi_{1,18},\phi_{2,4}),$
$[\phi_{2,1},\phi_{2,13}],$
$[\phi_{3,6}].$

$q=\zeta_{5}$  : $[\phi_{1,0},\phi_{3,8}],$
$[\phi_{2,4},\phi_{1,18}].$

$q=\zeta_{7}$  : $[\phi_{1,0},\phi_{2,10}].$

$q=\zeta_{8}$  : $[\phi_{1,0}],$
 $[\phi_{2,1}],$
 $[\phi_{2,4}],$
$[\phi_{3,2}].$

$q=\zeta_{9}$  : $[\phi_{1,0},\phi_{1,12}].$

$q=\zeta_{12}$  : $[\phi_{1,0},\phi_{1,18}],$
$[\phi_{2,1},\phi_{2,13}],$
$[\phi_{2,4},\phi_{2,10}].$\\
$ $\\
(vi) Let
$$\Hecke_q(G_8)=\left\langle s,t\,\left|\, sts=tst, 
\begin{array}{l}
(s-1)(s+1)(s+q)(s-q^5)=0,\\ (t-1)(t+1)(t+  q)(t-q^5)=0\end{array} \right.\right\rangle.$$
In this case, all Schur elements divide $s_{\phi_{1,18}}$.
The characters which, together with all characters of defect $0$, form an optimal basic set for $\Hecke_q(G_8)$ are:

$q= \pm 1$  : $(\phi_{1,0},\phi_{1,6},\phi_{2,1}).$

$q=\zeta_{3}$  : $(\phi_{1,0},\phi_{1,6},\phi_{1,12},\phi_{1,18}),$
$[\phi_{2,7''},\phi_{2,7'}].$

$q=\zeta_{5}$  : $\{\phi_{1,0},\phi_{2,1}\},$
$[\phi_{2,4}],$
 $[\phi_{3,2}].$

$q=\zeta_{6}$  : $(\phi_{1,0},\phi_{1,6},\phi_{1,12},\phi_{1,18}),$
$[\phi_{2,4},\phi_{2,10}].$

$q=\zeta_{7}$  : $[\phi_{3,2},\phi_{1,18}].$

$q=\zeta_{8}$  : $[\phi_{1,12}],$
$[\phi_{2,4}],$
$[\phi_{2,7''}],$
$[\phi_{3,2}].$

$q=\zeta_{9}$  : $[\phi_{2,4},\phi_{1,18}],$
$[\phi_{2,1},\phi_{2,13}].$

$q=\zeta_{10}$  : $\{\phi_{1,6},\phi_{2,1}\},$
$[\phi_{3,2}],$
$[\phi_{2,7''}].$

$q=\zeta_{12}$  : $[\phi_{1,12},\phi_{1,18}],$
$[\phi_{2,4},\phi_{2,10}],$
$[\phi_{2,7'},\phi_{2,7''}].$

$q=\zeta_{14}$  : $[\phi_{3,2},\phi_{1,18}].$

$q=\zeta_{15}$  : $[\phi_{1,6},\phi_{1,18}].$

$q=\zeta_{18}$  : $[\phi_{2,7''},\phi_{1,18}],$
$[\phi_{2,1},\phi_{2,13}].$

$q=\zeta_{30}$  : $[\phi_{1,0},\phi_{1,18}].$\\
$ $\\
(vii) Let
$$\Hecke_q(G_8)=\left\langle s,t\,\left|\, sts=tst, 
\begin{array}{l}
(s-1)(s+q^4)(s-q^5)(s+q^5)=0,\\ (t-1)(t+q^4)(t-  q^5)(t+q^5)=0\end{array} \right.\right\rangle.$$
In this case, all Schur elements divide $s_{\phi_{1,0}}$.
The characters which, together with all characters of defect $0$, form an optimal basic set for $\Hecke_q(G_8)$ are:

$q=-1$  : $(\phi_{1,0},\phi_{1,6},\phi_{2,1}).$

$q=\zeta_{3}$  : $(\phi_{1,0},\phi_{1,6},\phi_{1,12},\phi_{1,18}),$
$[\phi_{2,4},\phi_{2,10}].$

$q=\zeta_{5}$  : $\{\phi_{1,0},\phi_{2,7'}\},$
 $[\phi_{2,1}],$
 $[\phi_{3,4}].$

$q=\zeta_{7}$  : $[\phi_{1,0},\phi_{3,8}].$

$q=\zeta_{8}$  : $[\phi_{1,0}],$
$[\phi_{2,4}],$
$[\phi_{2,7'}],$
$[\phi_{3,6}].$

$q=\zeta_{9}$  : $[\phi_{1,0},\phi_{2,7''}],$
$[\phi_{2,1},\phi_{2,13}].$

$q=\zeta_{12}$  : $[\phi_{1,0},\phi_{1,6}],$
$[\phi_{2,4},\phi_{2,10}],$
$[\phi_{2,7'},\phi_{2,7''}].$

$q=\zeta_{15}$  : $[\phi_{1,0},\phi_{1,18}].$\\
$ $\\
(viii) Let
$$\Hecke_q(G_8)=\left\langle s,t\,\left|\, sts=tst, 
\begin{array}{l}
(s-1)(s+q)(s+q^4)(s-q^5)=0,\\ (t-1)(t+q)(t+  q^4)(t-q^5)=0\end{array} \right.\right\rangle.$$
The characters which, together with all characters of defect $0$, form an optimal basic set for $\Hecke_q(G_8)$ are:

$q=\pm 1$  : $(\phi_{1,0},\phi_{1,6},\phi_{2,1}).$

$q=\zeta_{3}$  : $(\phi_{1,0},\phi_{1,6},\phi_{1,18}).$

$q=\zeta_{5}$  : $\{\phi_{1,0},\phi_{3,2}\},$
$[\phi_{2,1}],$
$[\phi_{2,4}].$

$q=\zeta_{6}$  : $(\phi_{1,0},\phi_{1,6},\phi_{1,12},\phi_{1,18}),$
$[\phi_{2,1},\phi_{2,13}].$

$q=\zeta_{7}:$ $[\phi_{1,6}, \phi_{2,13}]$,   $[\phi_{1,12}, \phi_{2,1}]$.

$q=\zeta_{8}$  : $[\phi_{1,0}],$
$[\phi_{1,6}],$
$[\phi_{3,2}],$
$[\phi_{3,4}],$
$(\phi_{2,1}).$

$q=\zeta_{9}$  : $[\phi_{1,0},\phi_{2,13}],$
$[\phi_{2,1},\phi_{1,18}].$

$q=\zeta_{10}$  : $[[\phi_{1,0},\phi_{2,7''},\phi_{1,18}]].$

$q=\zeta_{12}$  : $[\phi_{1,0},\phi_{1,12}],$
$[\phi_{2,1},\phi_{2,13}],$
$[\phi_{1,6},\phi_{1,18}].$

$q=\zeta_{18}:$ $[\phi_{1,6}, \phi_{1,12}]$.

$q=\zeta_{24}:$ $[\phi_{2,1}, \phi_{2,13}]$.

$q=\zeta_{30}$  : $[\phi_{1,0},\phi_{1,18}].$

\subsubsection{Optimal basic sets for the group $G_9$}
$ $\\
\\
 Let
$$\Hecke_q(G_9)=\left\langle s,t\,\left|\, ststst=tststs, 
\begin{array}{l}
(s-1)(s-q^4)=0,\\ (t-1)(t- q^2)(t-  q^4)(t-q^6)=0\end{array} \right.\right\rangle.$$
The characters which, together with all characters of defect $0$, form an optimal basic set for $\Hecke_q(G_9)$ are:

$q=-1$  : $(\phi_{1,0},\phi_{2,5}),$
$(\phi_{2,1}).$

$q=\zeta_{3}$  : $(\phi_{1,0},\phi_{1,6},\phi_{1,12'},\phi_{1,24}),$
 $\{\phi_{2,5},\phi_{2,4}\},$
$[\phi_{1,12''},\phi_{1,18'}],$
$[\phi_{3,2}],$
$[\phi_{3,6'}].$

$q=i$  : $(\phi_{1,0},\phi_{2,10}),$
$(\phi_{1,6},\phi_{2,4}),$
$(\phi_{2,1}),$
$(\phi_{2,5}).$

$q=\zeta_{5}$  : $[[\phi_{1,0},\phi_{2,11''},\phi_{1,30}]],$
$[\phi_{1,6},\phi_{2,8}],$
 $[\phi_{2,1},\phi_{1,24}].$

$q=\zeta_{7}$  : $[\phi_{1,0},\phi_{2,8}],$
$[\phi_{2,1},\phi_{1,30}].$

$q=\zeta_{9}$  : $[\phi_{1,0},\phi_{1,18'}],$
$[\phi_{1,6},\phi_{1,24}],$
$[\phi_{1,12''},\phi_{1,30}].$

$q=\zeta_{12}$  : $[[\phi_{2,10},\phi_{1,0},\phi_{1,24}]],$
$[[\phi_{1,6},\phi_{1,30},\phi_{2,4}]],$
$[\phi_{2,1},\phi_{2,8}],$
$[\phi_{2,5},\phi_{2,13}].$

$q=\zeta_{15}$  : $[\phi_{1,0},\phi_{1,30}].$

$q=\zeta_{20}$  : $[\phi_{1,0},\phi_{3,12}],$
$[\phi_{3,2},\phi_{1,30}].$

$q=\zeta_{24}$  : $[\phi_{1,0},\phi_{1,24}],$
$[\phi_{1,6},\phi_{1,30}],$
$[\phi_{2,1},\phi_{2,13}],$
$[\phi_{2,5},\phi_{2,8}].$

\subsubsection{Optimal basic sets for the group $G_{10}$}
$ $\\
\\
 Let
$$\Hecke_q(G_{10})=\left\langle s,t\,\left|\, stst=tsts, 
\begin{array}{l}
(s-1)(s+q^2)(s-q^4)=0,\\ (t-1)(t- q^3)(t+  q^3)(t-q^6)=0\end{array} \right.\right\rangle.$$
The characters which, together with all characters of defect $0$, form an optimal basic set for $\Hecke_q(G_{10})$ are:

$q=-1:$ $(\phi_{ 1, 0 }, \phi_{ 1, 6 }, \phi_{ 1, 8 }),\, (\phi_{ 2, 1 }, \phi_{ 2, 8 }),\,  
                  (\phi_{ 3, 2 }),\, (\phi_{ 3, 10'' })$.
                   
$q=\zeta_3:$ $(\phi_{ 1, 0 }, \phi_{ 1, 12 }, \phi_{ 2, 9 }),\, (\phi_{ 1, 8 }, \phi_{ 2, 8 }),\,        
                        (\phi_{ 2, 12 }),\,  (\phi_{ 2, 5 }),\, (\phi_{ 2, 4 }),\,  (\phi_{ 1, 16 }, \phi_{ 1, 28 }, \phi_{ 2, 1 })$.
  
$q=i:$ $(\phi_{ 1, 0 }, \phi_{ 1, 18 }),\, (\phi_{ 1, 6 }, \phi_{ 1, 12 }),\, 
                 (\phi_{ 2, 4 }, \phi_{ 2, 10 }),\, (\phi_{ 2, 1 }, \phi_{ 2, 13 })$.
                 
$q=\zeta_5:$ $[\phi_{ 1, 8 }, \phi_{ 1, 22 }],\,[\phi_{ 1, 6 }, \phi_{ 1, 26 }],\, 
                          [[\phi_{ 1, 0 }, \phi_{ 1, 14 }, \phi_{ 1, 34 }]]$.
                          
$q=\zeta_7:$ $[\phi_{ 1, 0 }, \phi_{ 1, 28 }],\, [\phi_{ 1, 12 }, \phi_{ 1, 34 }]$.

$q=\zeta_8:$ $[[ \phi_{ 2, 15'' }, \phi_{ 1, 0 }, \phi_{ 1, 26 }]],\, 
                           [[\phi_{ 1, 8 },  \phi_{ 1, 34 }, \phi_{ 2, 7'' }]],\,$ 
                           
                           $\,\,\,\,\,\,\,\,\,\,\,\,\,\,\,\,\,\,\,\,\,[\phi_{ 2, 8 }, \phi_{ 2, 18 }],\, [\phi_{ 2, 5 }, \phi_{ 2, 21 }],\, 
                           [\phi_{ 2, 4 }, \phi_{ 2, 14 }],\, [\phi_{ 2, 1 }, \phi_{ 2, 17 }]$.
                           
$q=\zeta_9:$ $[\phi_{ 1, 0 }, \phi_{ 3, 8''}],\, [\phi_{ 2, 8 }, \phi_{ 1, 26 }],\, 
                            [\phi_{ 1, 12 }, \phi_{ 2, 18 }],\,$ 
                            
                            $\,\,\,\,\,\,\,\,\,\,\,\,\,\,\,\,\,\,\,\,\,[\phi_{ 3, 10'' }, \phi_{ 1, 34 }],\, 
                             [\phi_{ 1, 8 }, \phi_{ 2, 17 }],\, [\phi_{ 2, 1 }, \phi_{ 1, 28 }]$. 
  
$q=\zeta_{15}:$ $[\phi_{ 1, 0 }, \phi_{ 2, 18 }],\, [\phi_{ 2, 1 }, \phi_{ 1, 34 }]$.

$q=\zeta_{20}:$ $[\phi_{ 1, 0 }, \phi_{ 1, 34 }],\,[\phi_{ 2, 1 }, \phi_{ 2, 21 }],\, [\phi_{ 2, 4 }, \phi_{ 2, 18 }]$. 
  
 $q=\zeta_{24}:$  $[\phi_{ 1, 0 }, \phi_{ 2, 15''}],\, [\phi_{ 2, 7''}, \phi_{ 1, 34 }],\, 
                                  [\phi_{ 1, 8 }, \phi_{ 3, 8''}],\, [\phi_{ 3, 10'' },\phi_{ 1, 26 }]$.

\subsubsection{Optimal basic sets for the group $G_{12}$}
$ $\\
\\
 Let
$$\Hecke_q(G_{12})=\left\langle s,t,u\,\left|\, stus=tust=ustu,
\begin{array}{l}
(s-1)(s-q^2)=0,\\ (t-1)(t- q^2)=0,\\ (u-1)(u-q^2)=0\end{array} \right.\right\rangle.$$
The characters which, together with all characters of defect $0$, form an optimal basic set for $\Hecke_q(G_{12})$ are:

$q=-1:$ $\langle \phi_{1,0}, \phi_{2,4} \rangle$.

$q=\zeta_8:$ $[[\phi_{1,0},\phi_{1,12}, \phi_{2,1}]]$.

$q=\zeta_{12}:$ $[\phi_{1,0},\phi_{1,12}]$.

$q=\zeta_{24}:$ $[\phi_{1,0},\phi_{1,12}]$.

\subsubsection{Optimal basic sets for the group $G_{16}$}
$ $\\
\\
Let
$$\Hecke_q(G_{16})=\left\langle s,t\,\left|\, sts=tst, 
\begin{array}{l}
(s-1)(s-q)(s-q^2)(s-q^3)(s-q^4)=0,\\ (t-1)(t-q)(t-q^2)(t-q^3)(t-q^4)=0\end{array} \right.\right\rangle.$$
The characters which, together with all characters of defect $0$, form an optimal basic set for $\Hecke_q(G_{16})$ are:

$q=1:$    $ (\phi_{1,0},   \phi_{3,2}),\,( \phi_{2,1}, \phi_{4,3})$.

$q=-1:$   $ (\phi_{1,0},   \phi_{1,12}, \phi_{2,1}),
(\phi_{2,7},   \phi_{4,8}),\, [\phi_{4,6}],\,[\phi_{4,11}]$.

$q=\zeta_3:$   $ (\phi_{1,0},   \phi_{1,12}, \phi_{3,2}),\,( \phi_{2,1}, \phi_{2,7}, \phi_{2,13''}),\,
[\phi_{2,13'},   \phi_{2,25'}],\,$ 

$\,\,\,\,\,\,\,\,\,\,\,\,\,\,\,\,\,\,\,\,\,[\phi_{3,6}],\, [\phi_{3,10'}],\, [\phi_{6,5}],\, [\phi_{6,7}]$.

$q=i:$   $ (\phi_{1,0},   \phi_{2,7}, \phi_{3,6}),\,[[ \phi_{1,12}, \phi_{2,19'}, \phi_{1,36}]],\,$

$\,\,\,\,\,\,\,\,\,\,\,\,\,\,\,\,\,\,\,\,\,\{\phi_{2,1},   \phi_{4,8}\},\, \{\phi_{3,2}, \phi_{3,10''}\},\, \{\phi_{2,13'},\phi_{4,3}\}$.

$q=\zeta_{6}:$  $ (\phi_{1,0},   \phi_{1,12}, \phi_{1,24}, \phi_{1,36}, \phi_{1,48}),\,
[[ \phi_{2,7}, \phi_{2,19''}, \phi_{2,31}]],\,$

$\,\,\,\,\,\,\,\,\,\,\,\,\,\,\,\,\,\,\,\,\,\,\,[\phi_{2,13'},   \phi_{2,25'}],\, [\phi_{3,6}, \phi_{3,18''}],\, [\phi_{3,10''},\phi_{3,22}]$.

$q=\zeta_7:$  $ [\phi_{1,0},  \phi_{3,18''}],\, [\phi_{3,10''}, \phi_{1,48}],\,
[\phi_{2,1},   \phi_{4,17}],\, [\phi_{4,11}, \phi_{2,37}]$.

$q=\zeta_{8}:$  $ [[\phi_{1,0}, \phi_{2,19''},  \phi_{1,48}]],\, [\phi_{1,12}, \phi_{2,31}],\,
[\phi_{2,7}, \phi_{1,36}],$

$\,\,\,\,\,\,\,\,\,\,\,\,\,\,\,\,\,\,\,\,\,\,\,[\phi_{2,1},   \phi_{4,20}],\, [\phi_{3,2}, \phi_{3,18'}],\, [\phi_{3,10'},\phi_{3,26}]\,[\phi_{2,37},\phi_{4,8}]$.

$q=\zeta_9:$  $ [\phi_{1,0},  \phi_{3,26}],\, [\phi_{3,2}, \phi_{1,48}],\,
[\phi_{2,1},   \phi_{2,37}],\,[\phi_{2,7},   \phi_{2,25'}],\, [\phi_{2,13'}, \phi_{2,31}]$.

$q=\zeta_{10}:$   $ [[  \phi_{2,25'},\phi_{1,0}, \phi_{2,25''}]],\, [[\phi_{3,2}, \phi_{2,37}, \phi_{1,12}]],\,$

$\,\,\,\,\,\,\,\,\,\,\,\,\,\,\,\,\,\,\,\,\,\,\,[[  \phi_{1,36},\phi_{2,1}, \phi_{3,26}]],\, [[\phi_{2,13'},\phi_{1,48},\phi_{2,13''}]]$.

$q=\zeta_{12}:$  $ [[\phi_{1,0},  \phi_{1,24}, \phi_{1,48}]],\, [\phi_{1,12}, \phi_{1,36}],\,$

$\,\,\,\,\,\,\,\,\,\,\,\,\,\,\,\,\,\,\,\,\,\,\,[\phi_{2,1},   \phi_{2,25''}],\, [\phi_{2,13''}, \phi_{2,37}],\, [\phi_{3,2},\phi_{3,26}]$.

$q=\zeta_{14}:$  $ [\phi_{1,0},  \phi_{2,37}],\, [\phi_{2,1}, \phi_{1,48}]$.

$q=\zeta_{15}:$  $ [\phi_{1,0},  \phi_{4,20}],\, [\phi_{1,12}, \phi_{4,12}],\,
[\phi_{4,6},   \phi_{1,36}],\,$

$\,\,\,\,\,\,\,\,\,\,\,\,\,\,\,\,\,\,\,\,\,\,\,[\phi_{4,8}, \phi_{1,48}],\,
[\phi_{2,1},   \phi_{2,31}],\, [\phi_{2,7}, \phi_{2,37}]$.

$q=\zeta_{18}:$  $ [\phi_{1,0},  \phi_{1,36}],\, [\phi_{1,12}, \phi_{1,48}],\, [\phi_{2,1}, \phi_{2,37}]$.

$q=\zeta_{20}:$  $ [\phi_{1,0},  \phi_{5,8}],\, [\phi_{5,4}, \phi_{1,48}],\,
[\phi_{2,7},   \phi_{3,22}],\, [\phi_{3,6}, \phi_{2,31}]$.

$q=\zeta_{24}:$  $[\phi_{1,0},  \phi_{1,48}]$.

$q=\zeta_{30}:$  $ [\phi_{1,0},  \phi_{4,15}],\, [\phi_{4,3}, \phi_{1,48}],\,
[\phi_{2,1},   \phi_{3,26}],\, [\phi_{3,2}, \phi_{2,37}]$.

\subsubsection{Optimal basic sets for the group $G_{20}$}
$ $\\
\\
Let
$$\Hecke_q(G_{20})=\left\langle s,t\,\left|\, ststs=tstst, 
\begin{array}{l}
(s-1)(s-q)(s-q^2)=0,\\ (t-1)(t-q)(t-q^2)=0\end{array} \right.\right\rangle.$$
The characters which, together with all characters of defect $0$, form an optimal basic set for $\Hecke_q(G_{20})$ are:

$q=1:$  $ (\phi_{1,0}, \phi_{4,6}),\, (\phi_{2,1}, \phi_{2,7})$.

$q=-1:$   $ (\phi_{1,0},   \phi_{2,1}, \phi_{2,7}),\, [\phi_{4,3}],\, [\phi_{4,6}]$.

$q=i:$ $ [[\phi_{1,0}, \phi_{4,11}, \phi_{1,40} ]]$.

$q=\zeta_{5}:$  $[[\phi_{1,0}, \phi_{3,10'}, \phi_{1,40}]],\, [[\phi_{2,1},  \phi_{2,11}, \phi_{2,21}]],\,
[\phi_{2,7}, \phi_{2,27}]$.

$q=\zeta_{6}:$ $[[\phi_{2,21}, \phi_{1,0}, \phi_{2,27}]],\,  [[ \phi_{2,1}, \phi_{1,40}, \phi_{2,7}]]$.

$q=\zeta_{10}:$ $[[ \phi_{1,0},\phi_{1,20},  \phi_{1,40}]],\, [\phi_{2,1},  \phi_{2,21}]$.

$q=\zeta_{12}:$  $ [\phi_{1,0}, \phi_{5,12}],\,  [\phi_{5,4}, \phi_{1,40}]$.

$q=\zeta_{15}:$ $[\phi_{1,0},  \phi_{3,2}],\, [\phi_{2,7},\phi_{4,13}],\, 
[\phi_{4,3}, \phi_{2,27}],\, [\phi_{1,40}, \phi_{3,12}]$.

$q=\zeta_{20}:$  $ [\phi_{1,0},  \phi_{1,40}]$.

$q=\zeta_{30}:$  $[\phi_{1,0}, \phi_{2,21}],\, [\phi_{1,40}, \phi_{2,1}]$.

\subsubsection{Optimal basic sets for the group $G_{22}$}
$ $\\
\\
 Let
$$\Hecke_q(G_{22})=\left\langle s,t,u\,\left|\, stust=tustu=ustus,
\begin{array}{l}
(s-1)(s-q^2)=0,\\ (t-1)(t- q^2)=0,\\ (u-1)(u-q^2)=0\end{array} \right.\right\rangle.$$
The characters which, together with all characters of defect $0$, form an optimal basic set for $\Hecke_q(G_{22})$ are:

$q=-1:$    $ (\phi_{1,0},   \phi_{2,13}, \phi_{2,1})$ .

$q=\zeta_6:$ $ [[\phi_{1,0}, \phi_{4,6},  \phi_{1,30}]]$.

$q=\zeta_{10}:$  $ [[\phi_{1,0}, \phi_{2,13}, \phi_{1,30}]],\,
[\phi_{3,2}, \phi_{3,12}]$.

$q=\zeta_{12}:$ $ [\phi_{1,0}, \phi_{4,3}],\,  [\phi_{4,9},  \phi_{1,30}]$.

$q=\zeta_{20}:$  $[\phi_{1,0},  \phi_{3,12}],\, [\phi_{3,2}, \phi_{1,30}]$.

$q=\zeta_{30}:$  $ [\phi_{1,0},  \phi_{1,30}]$.

\subsection{Conjectures and remarks}

The following two conjectures have been 
checked to hold in all the cases presented in \S \ref{res}. Moreover, Conjecture \ref{gdraru}
has been 
proved for type $A$ and all exceptional Weyl groups in 
\cite{ChMy}, using  respectively \cite[3.7 \& 4.5]{Br} and the GAP3 package CHEVIE.
\footnote{For types $A$, $B$, etc.,
 M.~Fayers wrote to the second author a sketch of a proof
for Conjecture~\ref{gdraru} using purely combinatorics.}

Let $W$ be a complex reflection group.
Let $\theta: q \mapsto \xi$ be a specialization of $\Hecke_q(W)$ and let $(d_{\chi,\phi}^\theta)_{\chi \in \mathrm{Irr}(W), \phi \in \mathrm{Irr}(\Hecke_\xi)}$ be the corresponding decomposition matrix. 
\begin{conjecture}
For all $\phi \in \mathrm{Irr}(\Hecke_\xi)$, we have
$$\theta \left(\sum_{\chi \in \mathrm{Irr}(W)}
\frac{d_{\chi,\phi}^\theta}{s_\chi} \right) \neq 0$$
i.e., with the notation of Lemma \emph{\ref{GR}}, $\sum_{\chi \in \mathrm{Irr}(W)}
{d_{\chi,\phi}^\theta D_\chi^P}$ has as many factors of the form $(q-\xi)$ as $P(y)$.
\end{conjecture}

\begin{conjecture}\label{gdraru}
If $\chi, \psi \in \mathrm{Irr}(W)$ belong to the same block,
then
\begin{center}
$\theta \left( {s_\chi}/{s_\psi} \right) \in \mathbb{R}^\times.$
{\footnote{In fact, in all examples we have, the ratio belongs to $\mathbb{Q}^\times$.}}
\end{center}
\end{conjecture}

We would like to call the above ratio $\theta(s_\chi/s_{\psi})$
a {\emph{Brou\'e invariant}}.
The reason why we call it this is that in \cite[3.7 \& 4.5]{Br} 
M.~Brou\'e has shown that these ratios are Morita (respectively derived)
invariants (respectively up to sign) for block algebras with a suitable modular system.
\footnote{In \cite[3.7 \& 4.5]{Br}, M. Brou\'e considered the case of
 group algebras, but his argument is also valid for Hecke algebras.
Here, the ratio $\theta(s_{\chi}/s_{\psi})$ is slightly different from
 Brou\'e's original definition, but the interpretation between his
 and ours is clear if we treat Brou\'e's original ratios in \cite{Br} as 
specialized Schur element ratios and Conjecture~\ref{gdraru} is true.
The second author would like to thank M.~Brou\'e
 for explaining the proof of \cite[3.7 \& 4.5]{Br} to him.}

The authors expect that
 the parity condition on $\Ext^1$ for an appropriate quasihereditary
 cover,
 namely ${\mathcal{O}}={\mathcal{O}}_{\bf c}(W,V)$ of the rational Cherednik algebra
(see \cite{GGOR} for the definition of the rational Cherednik algebra and ${\mathcal{O}}$), holds.
Here, $V$ is the natural representation of $W$ and ${\bf c}$ is a charge for
the rational Cherednik algebra which is consistent with
 the parameters of the Hecke algebra $\Hecke_\xi$ via the ${\mathrm{KZ}}$-functor.

More precisely, 
let ${\mathcal{B}}$ be a block of $\Irr (W)$, which is defined with respect to the specialization
 $\theta: q \mapsto \xi$. Let $\chi_{\mathcal{B}}$ be the element in ${\mathcal{B}}$ with minimal
 $a$-value, which is conjecturally unique.
{\footnote{As in \cite{Gsurvey}, if $W$ is real, this uniqueness is true.}}
Conjecture~\ref{gdraru} allows us to define a partition  ${\mathcal{B}}={\mathcal{B}}_+
 \sqcup {\mathcal{B}}_-$, where
\[
{\mathcal{B}}_+:=\left\{
\chi \in {\mathcal{B}}\,|\, \theta(s_{\chi_{\mathcal{B}}}/s_{\chi})>0
\right\}\,\,\,\text{and}\,\,\, 
{\mathcal{B}}_-:=\left\{
\chi \in {\mathcal{B}}\,|\, \theta(s_{\chi_{\mathcal{B}}}/s_{\chi})<0
\right\} .
\]
We denote by $L(\lambda)$ the simple module in ${\mathcal{O}}$
 associated with $\lambda \in \Irr (W)$.
Define the {\emph{length function modulo $2$}}, 
denoted by $\ell_{{\mathrm{mod}}~2}$,  on $\Irr W$ by
\[
\ell_{{\mathrm{mod}}~2}(\lambda) :=
\left\{
\begin{array}{cl}
 0 \in \FF_2& {\mbox{ if }} \lambda \in {\mathcal{B}}_+,\\
1 \in \FF_2& {\mbox{ if }} \lambda \in {\mathcal{B}}_-,\\
\end{array}
\right.
\]
for some block ${\mathcal{B}} \ni \lambda$.
{\footnote{For type $A$, this length function is consistent with Cline-Parshall-Scott
 length function in Lusztig's conjecture as well as in the abstract
 Kazhdan-Lusztig theory.
So, we expect that the length function we define will be useful for 
the existence of the strong Kazhdan-Lusztig theory for category
 ${\mathcal{O}}$ of rational Cherednik algebras in the sense of Cline-Parshall-Scott.
 }}
We can state the {\emph{Parity Conjecture}},
 that is,
for any $\lambda, \mu \in \Irr (W)$, we have
$
\Ext^1_{\mathcal{O}}(L(\lambda),L(\mu)) \neq \{0\}$  {{ only if }}
\begin{center}
$\ell_{{\mathrm{mod}}~2}(\lambda)\neq \ell_{{\mathrm{mod}}~2}(\mu)
{\mbox{ and }} \lambda {\mbox{ and }} \mu {\mbox{ are in the same block.}}$
{\footnote{The validity of Lusztig's conjecture for
 quantum general linear groups at roots of unity and the equivalences between module
 categories of  quantized Schur algebras
and the category ${\mathcal{O}}$ of rational Cherednik algebras of type $A$
yield that the Parity Conjecture is
 true for type $A$.}}
\end{center}

Finally, we want to make the following remark on the lifting argument:
 Let $p$ be a prime number. The classical Fong-Swan theorem says that
if a group $G$ is $p$-solvable then
all the simple ${{\mathbb{F}}_{p^k}} [G]$-modules are lifted for some $k>0$
(cf.~\cite[7.5]{NaTsu}).
Note that all complex reflection groups that we deal with in this paper, except for the groups $G_{16}$, $G_{20}$ and $G_{22}$,
are solvable.
So if
 the modular system with parameters is well-defined for
our specialized Hecke algebra $\Hecke_\xi$ in characteristic zero and
a group algebra ${\mathbb{F}}_{p^k}[W]$ in positive characteristic, then
the modular simple
module over the group algebra ${\mathbb{F}}_{p^k}[W]$ can be lifted to
       $\Hecke_\xi$, hence over the generic Hecke algebra by the
       commutativity of two different modular systems
(see \cite[p.~141]{AM} for the modular system with parameters).
In this paper, we have stronger results. 
The statement above does not say that all simple
       $\Hecke_\xi$-modules
are liftable. It only says that some of the simple
       $\Hecke_\xi$-modules are liftable.
However, as one of our main results, 
it turns that all simple $\Hecke_\xi$-modules for our choices
       of parameters are liftable.

\section*{Appendix}$ $\\
Let us work out in detail the example of $G_{12}$, when
$$\Hecke_q(G_{12})=\left\langle s,t,u\,\left|\, stus=tust=ustu,
\begin{array}{l}
(s-1)(s-q^2)=0,\\ (t-1)(t- q^2)=0,\\ (u-1)(u-q^2)=0\end{array} \right.\right\rangle$$
and $q=\zeta_8$. We will use the programme GAP3 and we need the packages
CHEVIE and VKCURVE. Everything can be found on Jean Michel's webpage:
\begin{center}
{\sf{
http://people.math.jussieu.fr/\~{}jmichel 
}}
\end{center}

We type:
{ \small \begin{verbatim}

gap> RequirePackage("chevie");
gap> RequirePackage("vkcurve");
gap> W:=ComplexReflectionGroup(12);;
gap> CharNames(W);
[ "phi{1,0}", "phi{1,12}", "phi{2,1}", "phi{2,4}", 
  "phi{2,5}", "phi{3,2}", "phi{3,6}", "phi{4,3}" ]
  
\end{verbatim} } 
  
Following the notation of characters in CHEVIE, each irreducible character is denoted by $\phi_{d,b}$, where $d$ is the dimension of the representation and $b$ is the valuation of its fake degree. Therefore, the irreducible characters $\phi_{1,0}$ and $\phi_{1,12}$ are $1$-dimensional, $\phi_{2,1}$, $\phi_{2,4}$ and $\phi_{2,5}$ are $2$-dimensional,
$\phi_{3,2}$ and $\phi_{3,6}$ are $3$-dimensional, and $\phi_{4,3}$ is $4$-dimensional.
We now define the algebra $\Hecke_q(G_{12})$ and check the values of the Schur elements for $q=\zeta_8$:
 
{ \small \begin{verbatim}
 
gap> q:=Mvp("q");; 
gap> H:=Hecke(W,[[1,q^2]]);;
gap> s:=SchurElements(H);;
gap> List(s,i->Value(i,["q",E(8)]));
[ 0, 0, 0, 4, -288, 0, 0, 0 ]

\end{verbatim}}

Due to Criterion (1) of Section \ref{criteria}, the characters  $\phi_{2,4}$ and $\phi_{2,5}$ are of defect $0$. Now we will see what happens with the remaining irreducible representations of $\mathcal{H}_q(G_{12})$ when $q=\zeta_8$:

{ \small \begin{verbatim}
 
gap>  R:=Representations(H);;  
gap>  r:=List(R,i->List(i,j->List(j,k->List(k,l->Value(l,["q",E(8)])))));;
gap> r[1];
[ [ [ 1 ] ], [ [ 1 ] ], [ [ 1 ] ] ]
gap> r[2];
[ [ [ I ] ], [ [ I ] ], [ [ I ] ] ]

\end{verbatim}}

Hence, the characters $\phi_{1,0}$ and $\phi_{1,12}$ correspond to distinct irreducible representations of $\mathcal{H}_{\zeta_8}$. Now we will check whether the $2$-dimensional representation $\phi_{2,1}$ has either $\phi_{1,0}$ or $\phi_{1,12}$ as subrepresentation. 
For this, it is enough to check
 whether the generating matrices 
$\mathrm{r}[3][1]$ and $\mathrm{r}[3][2]$ have a common eigenvector for either the eigenvalue $1$ or $i$. We have created the following (not very sophisticated) functions:

{ \small \begin{verbatim}

Eigen:=function(a,Y) local d,res,values,vectors,l,j,n;
 d:=DimensionsMat(a)[1];
 res:=rec(values:=[],vectors:=[]);
 for l in Y do
 if DeterminantMat(a-l*IdentityMat(d))=0 then
       n:=NullspaceMat(TransposedMat(a-l*IdentityMat(d)));
       Add(res.vectors,n);
       for j in [1..Length(n)] do Add(res.values,l);od;fi; od;
 res.vectors:=Concatenation(res.vectors);
 return res;
end;
 
FindCommonEigenvector:=function(a,b,Y) local Ea,Eb,res,p,v,l,i;
 Ea:=Eigen(a,Y); 
 Eb:=Eigen(b,Y);
 res:=[];
 for i in [1..Length(Ea.vectors)] do
        for l in Eb.values do
              if b*Ea.vectors[i]=l*Ea.vectors[i] then
               Add(res,[Ea.values[i],l,Ea.vectors[i]]);fi;od;od; 
 for i in [1..Length(Eb.vectors)] do
        for l in Ea.values do
              if a*Eb.vectors[i]=l*Eb.vectors[i] then
               Add(res,[l,Eb.values[i],Eb.vectors[i]]);fi;od;od;
 return Set(res);
end; 
\end{verbatim}}
Now, let $a,b$ be two square matrices and $Y$ a list of complex numbers. The function
$Eigen(a,Y)$ returns a record consisting of the elements of $Y$ which are eigenvalues of  $a$, and a list of corresponding eigenvectors for $a$. The function  $FindCommonEigenvector(a,b,Y)$ returns a list of common eigenvectors for the matrices $a$ and $b$ with corresponding eigenvalues in $Y$ (it is not necessary that these vectors correspond to the same eigenvalue for $a$ and $b$). Thus, we check whether $\mathrm{r}[3][1]$ and $\mathrm{r}[3][2]$ have a common eigenvector for either the eigenvalue $1$ or $i$ in the following way:

{ \small \begin{verbatim}

gap> FindCommonEigenvector(r[3][1],r[3][2],[1,E(4)]);               
[  ]
\end{verbatim}}
$ $\\
They don't, so $\phi_{2,1}$ is also irreducible. We shall now apply Lemma \ref{GR} on generic degrees, recalling that, in this case, the Poincar\'e polynomial $P$ is equal to the Schur element $\mathrm{s}[1]$.

{ \small \begin{verbatim}

gap> D:=List(s,i->s[1]/i);;
gap> List(D,i->Value(i,["q",E(8)]));
[ 1, 1, 4, 0, 0, -3, -3, 2 ]

\end{verbatim}}

Let $(d_{\chi,\psi}^\theta)_{\chi \in \mathrm{Irr}(G_{12}),\psi \in \mathrm{Irr}({\mathcal{H}_{\zeta_8}})}$ be the decomposition matrix of $\Hecke_q(G_{12})$ with respect to the specialization $\theta: q \mapsto \zeta_8$. In order to have
$$\sum_{\chi \in \mathrm{Irr}(G_{12})}
{d_{\chi,\psi}^\theta D_\chi^P}(\zeta_8)=0,$$ the decomposition matrix must look as follows:

$$\begin{array}{l|ccccc}
\phi_{1,0} &1& & & &\\
\phi_{3,2} &a&1& b& &\\
\phi_{4,3} &1&1&1& &\\
\phi_{2,1} & &1& & &\\
\phi_{3,6} &c &1 &d& &\\
\phi_{1,12} & & &1 & &\\
\phi_{2,4} & & & & 1 &\\
\phi_{2,5} & & & &  &1\\
\end{array}$$ 
\\
where $\{(a,b),\,(c,d)\}=\{(1,0),\,(0,1)\}$.
To find which is  which, we use the fact that the modular reductions of the irreducible characters of $\mathcal{H}_q(G_{12})$ can be written uniquely as $\mathbb{N}$-linear combinations of the
 irreducible characters of $\mathcal{H}_{\zeta_8}$. We type:

{ \small \begin{verbatim}

gap> T:=CharTable(H).irreducibles;;
gap> t:=List(T,i->List(i,j->Value(j,["q",E(8)])));;            
gap> t[6]=t[3]+t[1];
true
gap> t[7]=t[3]+t[2];
true

\end{verbatim}}

Thus, $(a,b)=(1,0)$ and $(c,d)=(0,1)$,
and the decomposition matrix is:

$$\begin{array}{l|ccccc}
\phi_{1,0} &1& & & &\\
\phi_{3,2} &1&1& & &\\
\phi_{4,3} &1&1&1& &\\
\phi_{2,1} & &1& & &\\
\phi_{3,6} & &1 &1 & &\\
\phi_{1,12} & & &1 & &\\
\phi_{2,4} & & & & 1 &\\
\phi_{2,5} & & & &  &1\\
\end{array}$$

$ $\\
$ $\\
{\bf Acknowledgements.}\,\,\,
The authors would like to thank M.~Brou\'e, G.~Malle and J.~Michel for
 fruitful discussions on the topics of this paper. In particular, they are indebted to
  G.~Malle and J.~Michel for providing them preliminary versions of their article \cite{MaMi}.
The second author would like to thank I.~Gordon for his invitation to
 University of Edinburgh.
This research was partially supported by the Japanese Ministry of Education,
 Science, Sports and Culture, Grant-in-Aid for Young Scientists (B), 19740011, 2008,
 and by the EPSRC grant EP/G04984X/1.

\end{document}